\documentclass[10pt, reqno]{amsart}
\usepackage[right=1.5cm, left=1.5cm, top=1.5cm, bottom=1.5cm]{geometry}
\usepackage[english, activeacute]{babel}
\usepackage{pifont}
\usepackage{amsmath,amsthm,amsxtra}
\usepackage{epsfig}
\usepackage{amssymb}
\usepackage{multicol}
\usepackage{latexsym}
\usepackage{amsfonts}

\usepackage{hyperref}
\usepackage{xparse}
\usepackage{tabularx} 
\usepackage{multicol}
\usepackage{color}
\usepackage{hyperref}
\usepackage{float}
\usepackage{graphicx}
\usepackage{subcaption}
\usepackage{listings}

\newcommand{\R}{\mathbb{R}}

\newcommand{\Com}{\mathbb{C}}

\newcommand{\al}{\alpha}
\newcommand{\bt}{\beta}
\newcommand{\ga}{\gamma}

\newcommand{\ve}{\varepsilon}

\newcommand{\sech}{\operatorname{sech}}

\newcommand{\re}{\operatorname{Re}}
\newcommand{\ima}{\operatorname{Im}}

\newtheorem{thm}{Theorem}[section]

\theoremstyle{remark}
\newtheorem{rem}{Remark}[section]

\newcommand{\be}{\begin{equation}}
\newcommand{\ee}{\end{equation}}
\newcommand{\bp}{\begin{proof}}
\newcommand{\ep}{\end{proof}}
\newcommand{\bel}{\begin{equation}\label}
\newcommand{\eeq}{\end{equation}}
\newcommand{\bea}{\begin{eqnarray}}
\newcommand{\eea}{\end{eqnarray}}
\newcommand{\bee}{\begin{eqnarray*}}
\newcommand{\eee}{\end{eqnarray*}}
\newcommand{\ben}{\begin{enumerate}}
\newcommand{\een}{\end{enumerate}}

\date{\today}

\usepackage[english]{babel}

\title[NLS-like breathers]{Review on the stability of the Peregrine and related breathers}

\author{Miguel A. Alejo}
\thanks{(M.A.A.) Departamento de Matem\'aticas. Universidad de C\'ordoba, C\'ordoba, Spain. e-mail: {\tt malejo@uco.es}.}

\author{Luca Fanelli}
\thanks{(L.F.) Dipartimento di Matematica "Guido Castelnuovo", Universit$\grave{\text{a}}$ di Roma - {\it La Sapienza}, 
e-mail: {\tt fanelli@mat.uniroma1.it}} 

\author{Claudio Mu\~noz}
\thanks{(C.M.) CNRS and Departamento de Ingenier\'ia Matem\'atica DIM,  Universidad de Chile,
e-mail: {\tt cmunoz@dim.uchile.cl}. Partially funded by CMM Conicyt PIA AFB170001 and Fondecyt 1150202} 

\begin{document}

\begin{abstract}
In this note, we review stability properties in energy spaces of three important nonlinear Schr\"odinger breathers: Peregrine,
Kuznetsov-Ma, and Akhmediev. More precisely, we show that these breathers are \emph{unstable} according to a standard definition
of stability. Suitable Lyapunov functionals are described, as well as their underlying spectral properties. As an immediate 
consequence of the first variation of these functionals, we also present the corresponding nonlinear ODEs fulfilled by these NLS breathers. The notion of 
global stability for each breather above mentioned is finally discussed. Some open questions are also briefly mentioned.
\end{abstract}

\maketitle
\numberwithin{equation}{section}

\vspace{-1cm}

\section{Introduction}

\begin{figure}[h!]
          \includegraphics[scale=0.24]{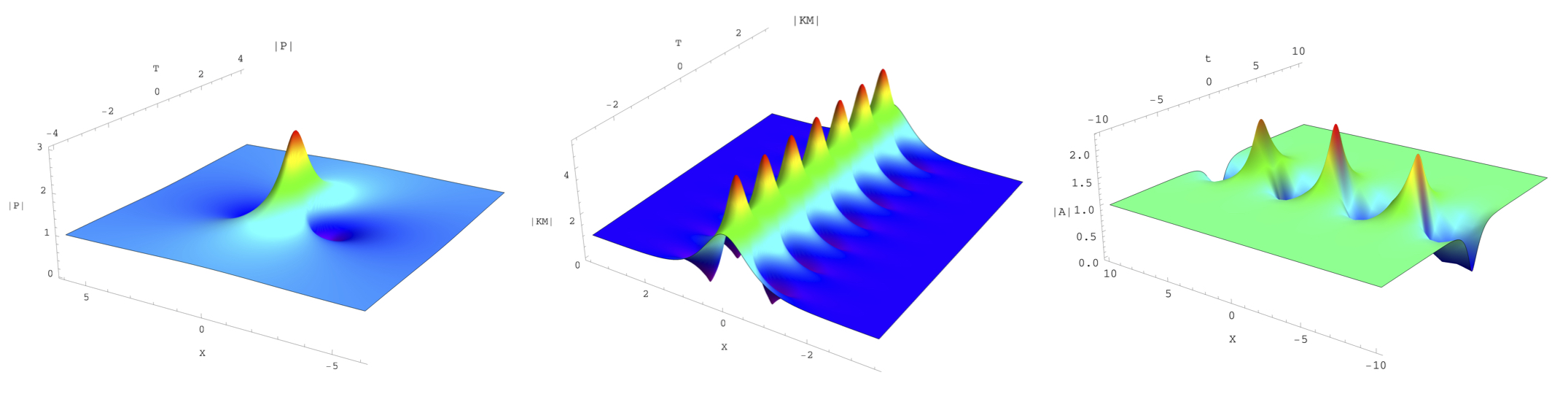}
              \caption{\small \emph{Left:} Absolute value $|B_P|$ of the Peregrine breather \eqref{P}. \emph{Center:} $|B_{KM}|$, \eqref{KM}. \emph{Right:} $|B_A|$ \eqref{Ak}.}\label{P_KM}
          \end{figure}

In this short review we describe a series of mathematical results related to the stability of the Peregrine breather and other explicit solutions to 
the cubic Nonlinear Schr\"odinger (NLS) equation, an important candidate to modelize rogue waves. The mentioned model is NLS posed on the real line
\be\label{NLS0}
i \partial_t u +\partial_x^2 u +|u|^2u=0, \quad u(t,x)\in\Com, \quad (t,x)\in \R^2. 
\ee
We assume a nonzero boundary value condition (BC) at infinity, in the form of an \emph{Stokes wave} $e^{it}$: for all $t\in\R$, 
\be\label{Stokes}
|u(t,x)-e^{it}| \to 0 \quad \hbox{as} \quad x\to \pm\infty.
\ee
It is well known that \eqref{NLS0} possesses a huge family of complex solutions. Among them, a fundamental role for the dynamics is played by {\it breathers}.  
We shall say that a particular smooth solution to \eqref{NLS0}-\eqref{Stokes} is a breather if, up to the invariances of the equation, 
its dynamics shows the evolution of some concentrated quantity in an oscillatory fashion. NLS has scaling, shifts, phase and Galilean invariances: namely, if $u$ solves \eqref{NLS0}, another solution to \eqref{NLS0} is
\be\label{standing_wave}
u_{c,v,\ga,x_0,t_0}(t,x):= \sqrt{c} \, u\left( c (t-t_0),\sqrt{c} (x-vt-x_0) \right) \exp \Big(  ict + \frac i2 xv  -\frac i4 v^2 t + i\ga \Big),
\ee
In this paper we review the known results about  stability in Sobolev spaces of the Peregrine ($P$) breather\footnote{Or Peregrine soliton, but because of the nature of its variational formulation, it is more a breather than a soliton.} \cite{Peregrine}
\be\label{P}
B_P(t,x):= 
e^{it}\Big(1-\frac{4(1+2it)}{1+ 4t^2 +2x^2}\Big);
\ee
we will also present, with less detail, the analogous properties of the Kuznetsov-Ma  ($KM$) and Akhmediev ($A$) breathers: 
(i) if $a>\frac12$, the Kuznetsov-Ma  (KM)  breather \cite{Kuznetsov,Ma,Akhmediev} (see Fig. \ref{P_KM} for details)
\be\label{KM}
\begin{aligned}
B_{KM}(t,x) : =  
e^{it}\Big(  1- \sqrt{2}\beta \frac{(\beta^2 \cos(\al t)  + i\al \sin(\al t)) }{ \al \cosh(\beta x) - \sqrt{2} \beta \cos(\al t)}\Big),\qquad  \al :=  (8a(2a-1))^{1/2}, \quad \beta := (2(2a-1))^{1/2},
\end{aligned}
\ee
and (ii) for $a\in (0,\frac12)$, the Akhmediev breather \cite{Akhmediev}
\be\label{Ak}
\begin{aligned}
B_A(t,x):=  e^{it}\Big(  1+ \frac{\alpha^2 \cosh( \beta t) +i\beta \sinh(\beta t) }{ \sqrt{2a} \cos(\alpha x) -\cosh(\beta t)}\Big),\qquad 
\beta:= (8a(1-2a))^{1/2}, \quad \alpha:=(2(1-2a))^{1/2}.
\end{aligned}
\ee

Notice the oscillating character of the three above examples. 
In addition, also notice that $B_{KM}$ is time-periodic, while $B_A$ is space-periodic. Although the Peregrine breather is not periodic in time, it is a particular limiting ``degenerate'' case of the two last cases.

NLS \eqref{NLS0} with nonzero BC \eqref{Stokes} is believed to describe the emergence of rogue or freak waves 
in deep sea \cite{Peregrine,Akhmediev2,Shrira}. Peregrine waves were experimentally observed 10 years ago in \cite{Kliber0}. The model itself is also a well-known example of the mechanism known as modulational instability 
\cite{Zakharov,Akhmediev2}. For an alternative explanation to the rogue wave phenomenon which is stable under perturbations, see \cite{Bona_Saut,BPSS}. 


Along these lines, we will explain that Peregrine and two other breathers are unstable according to a standard definition of stability. It could be the case that a less demanding definition of stability, involving infinite energy solutions, could repair this problem. However, such a question is still an open problem.

This paper is organized as follows. In Section \ref{Sect:2} we recall some standard results for NLS with zero and nonzero background,  and the notion of modulational instability and local well-posedness. Section \ref{Sect:3} is devoted to the conservation laws of NLS, and Section \ref{Sect:4} to the notion of stability. Sections \ref{Sect:5} and \ref{Sect:6} review the Peregrine, Kuznetsov-Ma and Akhmediev breathers' stability properties. Finally, Section \ref{Sect:7} is devoted to a discussion, final comments and conclusions.

\section{Modulational instability}\label{Sect:2}

\subsection{A quick review of the literature}

Let us briefly review the main results involving equation \eqref{NLS0} in the zero and nonzero BC case. A much more complete description of the current literature can be found in the papers present in this volume, and in \cite{Cazenave,Dauxois,Yang}.


NLS \eqref{NLS0} is a well-known integrable model, see \cite{ZS}, and describes 
the propagation of pulses in nonlinear media and gravity waves in the ocean \cite{Dauxois}. The local and global well-posedness theory for NLS with zero BC at infinity was initiated by Ginibre and Velo \cite{GV}, 
see also Tsutsumi \cite{Tsutsumi} and Cazenave and Weissler \cite{CW}. Finally, see Cazenave \cite{Cazenave} for a complete account on the different NLS equations. 
One should have in mind that \eqref{NLS0} is globally well-posed in $L^2(\R)$, which has been proved by Tsutsumi in \cite{TSU} and ill-posed in 
$H^s(\R)$, $s<0$, as showed in \cite{KPV}, where the authors prove the lack of uniform continuity of the solution map. 


In the zero background case, one has standard solitons for \eqref{NLS0}:
\be\label{soliton}
Q(t,x):=\sqrt{c} \sech(\sqrt{c}(x-vt-x_0)) e^{i \left( ct + \frac12 xv -\frac14  v^2 t +\gamma_0 \right)}, \quad c>0,~ v,x_0,\ga_0 \in\R.
\ee
These are time periodic, spatially localized solutions of \eqref{NLS0}, and orbitally stable, see Cazenave-Lions \cite{CL}, Weinstein \cite{Weinstein}, 
and Grillakis-Shatah-Strauss \cite{GSS}. See also \cite{Kap,MManihp,MMT2} 
for the case of several solitons. 

\subsection{Some heuristics} NLS with nonzero boundary condition, represented in \eqref{NLS0}-\eqref{Stokes}, is characteristic of the \emph{modulational instability} phenomenon, which -roughly speaking- says that small perturbations of the exact Stokes solution $e^{it}$ are unstable and grow quickly. This unstable behavior leads to a nontrivial competition with the (focusing) nonlinearity, time at which the solution is apparently stabilized.


There are plenty of works in the literature dealing with this phenomenon, not only in the NLS case. Usually also called Benjamin-Feir mechanism \cite{Bridges}, the NLS case has been described in a series of papers \cite{Akhmediev,Zakharov,Achilleos,Brunetti,Khawaja}. See also references therein for more details on the physical literature. Here we present the standard, simple but formal explanation of this phenomenon, in terms of a frequency analysis of the linear solution.


To this aim, consider \emph{localized} perturbations of \eqref{NLS0}-\eqref{Stokes}, of the form
\be\label{deco}
u(t,x) = e^{it}(1+ w(t,x)), \quad w~ \hbox{ unknown}.
\ee
Notice that this ansatz is motivated by \eqref{P}.
Then \eqref{NLS0} becomes a modified NLS equation with a zeroth order term, which is real-valued, and has the wrong sign:
\be\label{mNLS}
i\partial_t w + \partial_x^2 w   + 2 \re w + 2|w|^2 + w^2 +|w|^2w  =  0.
\ee
The associated linearized equation for \eqref{mNLS} is just
\be\label{phi_eqn}
\partial_t^2 \phi+\partial_x^4 \phi + 2 \partial_x^2 \phi =0, \qquad \phi =\re w.
\ee
This problem has some instability issues, as a standard frequency analysis reveals: looking for a formal standing wave $\phi=e^{i(kx-\omega t)}$ solution to \eqref{phi_eqn}, one has $\omega(k) = \pm |k| \sqrt{k^2 -2}$, which shows that for small wave numbers ($|k|<\sqrt{2}$) the linear equation behaves in an ``elliptic'' fashion, and exponentially (in time) growing modes are present from small perturbations of the vacuum solution. A completely similar conclusion is obtained working in the Fourier variable. This singular behavior is not present the equation is defocusing, that is \eqref{NLS0} with nonlinearity $-|u|^{2} u$.


This phenomenon is similar to the one present in the \emph{bad Boussinesq} equation, see Kalantarov and Ladyzhenskaya  \cite{KL}. However, in the latter case the situation is even more complicated, since the linear equation is ill-posed for all large frequencies, unlike NLS \eqref{NLS0}-\eqref{Stokes}, which is only badly behaved at small frequencies.

\subsection{Local well-posedness}

The above heuristics could lead to think that the model \eqref{NLS0}-\eqref{Stokes} is not well-posed (in the Hadamerd sense \cite{Hadamard}) in standard Sobolev spaces (appealing to physical considerations, we will only consider solutions to these models with
\emph{finite energy}). Recall that for $s\in\R$, the vector space $H^s(\R;\Com)$ corresponds to the Hilbert space of complex-valued functions 
$u:\R\to\Com$, such that $\int (1+\xi^2)^{s} |\hat u(\xi)|^2 d\xi<+\infty$, endowed with the standard norm (with the hat $\hat{\cdot}$ we denote the Fourier transform). Also, $H^s_\sharp:= H^s_\sharp((0,\frac{2\pi}{a}))$ denotes the Sobolev space $H^s$ of $\frac{2\pi}{a}$-space periodic functions. In \cite{Munoz,AFM2} it was shown that, even if  there is no time decay for the linear dynamics due to the modulationally unstable regime, the equation is still locally well-posed.

%

\begin{thm}\label{MT}
Let $s>\frac12$ and $a>0$. The  NLS with nonzero background \eqref{mNLS} is locally well-posed in $H^s$, in the aperiodic case, and in $H^s_\sharp$ in the periodic case. 
\end{thm}


The main feature in the proof of Theorem \ref{MT} is the fact that, if we work in Sobolev spaces, in principle \emph{there is no $L^1-L^\infty$ decay estimates for the linear dynamics}. Moreover, one has exponential growth in time of the $L^2$ norm, and therefore no suitable Strichartz estimates seems to be available, unless one cuts off some bad frequencies. Consequently, Theorem \ref{MT} is based in the fact that we work in dimension one, and that for $s>\frac12$, we have the inclusion $H^s \hookrightarrow  L^\infty$. 
See \cite{Achilleos} for early results on the Cauchy problem for \eqref{NLS0} in the periodic case, at high regularity ($H^2$ global weak solutions).


Note that $P$ in \eqref{P} is always well-defined, and has essentially no loss of regularity, confirming in some sense the intuition and the conclusions in Theorem \ref{MT}. Also, note that using the symmetries of the equation, we have LWP for any solution of \eqref{NLS0} of the form $u(t,x) = u_{c,v,\ga}(t,x) +w(t,x)$, for $w\in H^s$, and $s>\frac12$, with $u_{c,v,\ga}$ defined in \eqref{standing_wave}.


Two interesting questions are still open: global existence vs. blow-up, and ill-posedness of the flow map for lower regularities. Since \eqref{NLS0} is integrable, it can be solved, at least formally, by inverse scattering methods. Biondini and Mantzavinos \cite{BM} showed the existence and long-time behavior of a global solution to \eqref{NLS0} in the {integrable} case, under certain exponential decay assumptions at infinity, and a \emph{no-soliton} spectral condition. In this paper we have decided to present results stated just in some energy space, with no need of extra decay condition.


\section{Conserved quantities}\label{Sect:3}

Being an integrable model, \eqref{NLS0}-\eqref{Stokes} possesses an infinite number of conserved quantities \cite{ZS}. Here we review the most important for the question of stability: mass, energy and momentum. For both for $KM$ and $P$, one has the mass, momentum and energy
\be\label{Mass_NLS}
M[u]: =\int (|u|^2 -1), \quad P[u]:= \ima\int  (\bar u-e^{-it}) u_x, \quad E[u]:=\int |u_x|^2 - \frac12 \int (|u|^2-1)^2,
\ee
and the Stokes wave + $H^2$ perturbations conserved energy:
\be\label{F_NLS}
F[u]:= \int \Big(|u_{xx}|^2 -3 (|u|^2-1)|u_x|^2 -\frac12((|u|^2)_x)^2  + \frac12 (|u|^2-1)^3 \Big). 
\ee
In \cite{Munoz2}, it was computed the mass, energy and momentum of the $P$ \eqref{P} and $KM$ \eqref{KM} breathers. Indeed, one has \cite{Munoz2,AFM1}
\[
M[B_P]= E[B_P]= P[B_{P}]=F[ B_P]= P[B_{KM}]=0, \qquad M[B_{KM}] = 4\beta, \quad E[B_{KM}] =-\frac 83\bt^3, \quad F[B_{KM}]=\frac45\bt^5.
\]
We conclude that $KM$ and $P$ breathers are zero speed solutions. Note instead that, under a suitable Galilean transformation, they must have nonzero momentum. Note also that $P$ has same energy and mass as the Stokes wave solution (the nonzero background), a property not satisfied by the standard soliton on zero background. Also, compare the mass and energy of the Kuznetsov-Ma breather with the ones obtained in \cite{AM} for the mKdV breather.

Assume now that $u=u(t,x)$ is a $\frac{2\pi}{a}$-periodic solution to \eqref{NLS0}. Two standard conserved quantities for \eqref{NLS0} in the periodic setting are mass and energy
\be\label{Mass_Energy_per}
M_A[u]: =\int_{0}^{\frac{2\pi}{a}} (|u|^2 -1), \qquad E_A[u]:= \int_{0}^{\frac{2\pi}{a}}  \Big( |u_x|^2 - \frac12   (|u|^2-1)^2 \Big).
\ee
A third one, appearing from the integrability of the equation, is given by \cite{AFM1}
\be\label{F_NLS_per}
F_A[u]:= \int_{0}^{\frac{2\pi}{a}} \Big(|u_{xx}|^2 -3 (|u|^2-1)|u_x|^2 -\frac12((|u|^2)_x)^2  + \frac12 (|u|^2-1)^3 \Big). 
\ee

\section{Orbital stability}\label{Sect:4}

%

From a physical and mathematical point of view, understanding the stability properties of 
candidates to rogue waves is of uttermost importance because not all the observed patterns bear the same qualitative and quantitave information.

Since the equation is locally well-posed and does have continuous in time solutions, it is possible to define a notion of orbital stability for the Peregrine, Kuznetsov-Ma and Akhmediev breathers. To study the stability properties of such waves is key to validate them as candidates for explaining rogue waves, see \cite{Calini}. First, we consider the aperiodic case.

Mathematically speaking, the notion of orbital stability is the one to have in mind. Fix $s>\frac12$, and $t_0\in \R$. We say that a particular globally defined solution $U=e^{it}(1+W)$ of \eqref{NLS0} is \emph{orbitally stable} in $H^s$ if there are constants $C_0,\ve_0>0$ such that, for any $0<\ve<\ve_0$, if $\| w_0 -  W(t_0)\|_{H^s} <  \ve$, then 
\be\label{Stability}
\sup_{t\in \R}\inf_{y\in\R}  \|w(t) - W(t,x-y) \|_{H^s} <C_0\, \ve.
\ee
Here $w(t)$ is the solution to the IVP \eqref{mNLS} with initial datum $w(t_0)=w_0$, constructed in Theorem \ref{MT}, and $x_0(t)$ can be assumed continuous because the IVP is well-posed in a continuous-in-time Sobolev space. 


Note that no phase correction is needed in \eqref{Stability}: equation \eqref{mNLS} is no longer $U(1)$ invariant, and any phase perturbation of a modulationally unstable solution $u(t)$ in \eqref{NLS0}, of the form $u(t) e^{i\ga}$, $\ga\in \R$, \emph{requires an infinite amount of energy}. The same applies for Galilean transformations. If \eqref{Stability} is not satisfied, we will say that $U$ is unstable. Note additionally that condition \eqref{Stability} requires $w$  globally defined, otherwise $U$ is trivially unstable, since $U$ is globally defined.


Recall that NLS solitons on a zero background \eqref{soliton} satisfy \eqref{Stability} (with an additional phase correction) for $s=1$, see e.g. \cite{CL,GSS,Weinstein}. Some breather solutions of canonical integrable equations such as mKdV and Sine-Gordon have been shown stable using Lyapunov functional techniques, see \cite{AM,AM1,AMP1,AMP2,AMP3}. See also \cite{Gong1,Gong2} for a rigorous treatment using IST, and \cite{Dauxois,Peli_Yang,Yang} for more results for other canonical models, 
and \cite{GP1, GP2} for the stability of periodic waves and kinks for the defocusing NLS. 
For several years, a proof of stability/instability of NLS breathers was open, due to the difficult character (no particular sign) of conservation laws.

Now we consider an adapted version of stability for dealing with the Akhmediev breather \eqref{Ak}. We must fix a particular spatial period, which for the latter case will be settled as $L=\frac{2\pi}a$, for some fixed $a\in (0,\frac12)$, see \eqref{Ak}.


By stability in this case, we mean the following. Fix $s>\frac12$, and $t_0\in \R$. We say that a particular $\frac{2\pi}{a}$-periodic globally defined solution $U= e^{it}(1+W)$ of \eqref{NLS0} is \emph{orbitally stable} in $H^s_\sharp(\frac{2\pi}{a})$ if there are constants $C_0,\ve_0>0$ such that, for any $0<\ve<\ve_0$, if $\| u_0 -  U(t_0)\|_{H^s_\sharp} <  \ve$, then
\be\label{Stability_per}
 \sup_{t\in\R} \inf_{y,s\in\R}  \|u(t) - e^{is}U(t,x-y) \|_{H^s_\sharp} <C_0\, \ve.
\ee
If \eqref{Stability_per} is not satisfied, we will say that $U$ is unstable. Note how in this case phase corrections are allowed. This is because they are finite energy perturbations in the periodic case. In other words, any change of the form $B_A(t,x)e^{i\ga}$, $\ga\in\R$ of the Akhmediev breather $B_A(t,x)$ is a finite energy perturbation. The remaining sections of this review are devoted to show that all breathers considered in the introduction are unstable according to the previously introduced definitions.

\section{The Peregrine breather}\label{Sect:5}

Recall the Peregrine breather introduced in \eqref{P}. Note that it is a polynomially decaying (in space and time) perturbation of the nonzero background given by the Stokes wave $e^{it}$. 
Using a simple argument coming from the modulational instability of the equation \eqref{NLS0}, in \cite{Munoz2} it was proved that $B_{P}$ is unstable with respect to perturbations in Sobolev spaces $H^s$, $s>\frac12$. Previously, Haragus and Klein \cite{KH} showed numerical instability of the Peregrine breather, giving a first hint of its unstable character. 

\begin{thm}\label{Insta_Per}
The Peregrine breather \eqref{P} is unstable under small $H^s$ perturbations, $s>\frac12$.
\end{thm}

The proof of this result uses the fact that Peregrine breathers are in some sense converging to the background final state (i.e. they are asymptotically stable) in the whole space norm $H^s(\R)$, a fact forbidden in Hamiltonian systems with conserved quantities and stable solitary waves. 


Theorem \ref{Insta_Per} is in contrast with other positive results involving breather solutions \cite{AM,AM1,MP}. In those cases, the involved equations (mKdV, Sine-Gordon) were globally well-posed in the energy space (and even in smaller subspaces), with uniform in time bounds. Several physical and computational studies on the Peregrine breather can be found in \cite{Dysthe,Brunetti} and references therein. A recent stability analysis was performed in \cite{Zweck} in the case of complex-valued Ginzburg-Landau models. The proof of Theorem \ref{Insta_Per} is in some sense a direct application of the notion of modulational instability together with an asymptotic stability property. 

%
          
\subsection{Sketch of proof of Theorem \ref{Insta_Per}}
This proof is not difficult, and it is based in the notion of  \emph{asymptotic stability}, namely the convergence at infinity of perturbations of the breather. Fix $s>\frac12$. Let us assume that the Peregrine breather $P$ in \eqref{P} is orbitally stable, as in \eqref{Stability}. Write
\be\label{Q_def}
P(t,x)= e^{it} (1+ Q(t,x)),\qquad Q(t,x):= ~ -\frac{4(1+2it)}{1+ 4t^2 +2x^2}.
\ee
Now consider, as a perturbation of the Peregrine breather, the Stokes wave \eqref{Stokes}.  One has \cite{Munoz2}
\[
\lim_{t\to +\infty} \|  e^{it} -P(t) \|_{H^s}= \lim_{t\to +\infty} \| Q(t) \|_{H^s} =0.
\]
Therefore, we have two modulationally unstable solutions to \eqref{NLS0} that converge to the same profile as $t\to +\infty$. 
This fact contradicts the orbital stability, since for $y=x_0(t)\in \R$ given in \eqref{Stability}, 
$0< c_0:=\|Q(0,x-x_0(0))\|_{H^s}$ is a fixed number, but if $t_0=T$ is taken large enough, $ \| Q(T) \|_{H^s}$ can 
be made arbitrarily small. This proves Theorem \ref{Insta_Per}.


Although Theorem \ref{Insta_Per} clarifies the stability/instability question for the Peregrine breather, other questions remain unsolved. 
Is the Peregrine breather stable under less restrictive assumptions on the perturbed data? A suitable energy space for the Peregrine breather could be
\[
\mathcal E:= \left\{  u\in L^\infty(\R), \quad |u|^2 - 1 \in L^2(\R), \quad u_x\in L^2(\R) \right\}, 
\]
endowed with the metric $d(u_1,u_2):= \|u_{1}-u_{2}\|_{L^\infty}  +  \|u_{1,x}-u_{2,x}\|_{L^2} + \| |u_1|^2 -|u_2|^2\|_{L^2}.$
This space is standard for the study of kink structures in Gross-Pitaevski \cite{Gallo2,Gerard}. However, even in this space the argument used in Theorem \ref{Insta_Per} works, giving instability as well. In other words, the asymptotic stability property \emph{in the whole space} is key for the instability.

\subsection{Variational characterization}  In the following lines, we discuss some improvements of the previous result. In particular, we discuss the variational characterization of the Peregrine breather. For an introduction to this problem in the setting of breathers, see e.g. \cite{Munoz}. In \cite{AFM1}, the authors quantified in some sense the instability of the Peregrine breather.

\begin{thm}[Variational characterization of Peregrine]\label{TH1_Per}
Let $B=B_{P}$ be any Peregrine breather. Then $B$ is a \emph{critical point} of a real-valued functional $F[u]$ \eqref{F_NLS}, in the sense that
\be\label{Hp}
F'[B](z)=0, \quad \hbox{for all } ~ z\in H^2(\R;\Com). 
\ee
Moreover, $B$ satisfies the nonlinear ODE
\be\label{Ec_P}
\begin{aligned}
B_{(4x)} + 3B_x^2 \bar B +(4 |B|^2-3) B_{xx}+ B^2 \bar B_{xx}  + 2  |B_x|^2 B + \frac32 (|B|^2-1)^2 B  =0.
\end{aligned}
\ee
\end{thm}

Theorem \ref{TH1_Per} reveals that Peregrine breathers are, in some sense, degenerate. More precisely, contrary to other breathers, the characterization of $P$ does not require the mass and the energy, respectively. The absence of these two quantities may be related to the fact that $M[B_P]=E[B_P]=0$, meaning a particular form of instability (recall that mass and energy are convex terms aiding to the stability of solitonic structures). We would like to further stress the fact that the variational characterization of the famous Peregrine breather is in $H^2$, since mass and energy are useless. 


The proof of Theorem \ref{TH1_Per} is simple, variational and follows previous ideas presented in \cite{AM} for the case of mKdV breathers, and \cite{AMP1} for the case of the Sine-Gordon breather (see also \cite{MP} for a recent improvement of this last result, based in \cite{AM1}).  The main differences are in the complex-valued nature of the involved breathers, and the nonlocal character of the $KM$ and $P$ breathers. 


The following result gives a precise expression for the lack of stability in Peregrine breathers. 
Recall that $\sigma_c(\mathcal L)$ stands for the continuum spectrum of a densely defined unbounded linear operator $\mathcal L$. Essentially, the continuous spectrum of the second derivative of the Lyapunov functional $F''$ stays below zero. 

\begin{thm}[Direction of instability of the Peregrine breather]\label{TH4}
Let $B=B_P$ be a Peregrine breather, critical point of the functional $F_P$ defined in \eqref{F_NLS}. Then the following is satisfied: let $z_0\in H^2$ be any sufficiently small perturbation, and $w=w(t):= e^{-it} \partial_x z_0\in H^1$. Then, as $t\to +\infty$,
\be\label{Insta_P}
\begin{aligned}
F''[B_P](z_0,z_0)=  \frac12 \int  (|w_{x}|^2 -|w|^2 - w^2 )(t) +O(\|z_0\|_{H^1}^3) +o_{t\to +\infty}(1).
\end{aligned}
\ee 
\end{thm}

From \eqref{Insta_P}, one can directly check that $F''[B_P](z_0,z_0)<0$ for a continuum of small $z_0$ and large times.  
The proof of Theorem \ref{TH4} is consequence of the following identity. For each $z\in H^2(\R)$, we have 
\be\label{decomp}
F[ B_P + z]=  F[ B_P] + \mathcal G[z] + \mathcal Q[z] + \mathcal N[z],
\ee
where $F[ B_P]=0$, $\mathcal G[z] = 2\re\int \bar{z} G[B_P]=0$, with $G[B_{P}]=\eqref{Ec_P}$, and $\mathcal Q[z]$ a quadratic functional of the form $\mathcal Q[z] := \re \int\bar{z}\mathcal{L}_{P}[z]dx$, where $ \mathcal{L}_{P}[z]$ is a matrix linear operator \cite{AFM1}.
Finally, assuming $\|z\|_{H^1}$ small enough, the term $|\mathcal N[z]|$ is of cubic order and small. 


\section{The Kuznetsov-Ma and Akhmediev breathers}\label{Sect:6}

Here we describe stability properties of the other two important breathers for NLS: the Kuznetsov-Ma  (KM)  breather \eqref{KM} (see Fig. \ref{P_KM} 
and \cite{Kuznetsov,Ma,Akhmediev}  for details) and  the Akhmediev breather \eqref{Ak} (see \cite{Akhmediev}).
%

\subsection{Kuznetsov-Ma} Most of the results obtained in the Peregrine case are also available for the Kuznetsov-Ma breather. We start by noticing that $B_{KM}$ is, unlike Peregrine, a Schwartz perturbation of the Stokes wave $e^{it}$, and therefore a smooth classical solution of \eqref{NLS0}.  
It has been also observed in optical fibre experiments, see Kibler et al. \cite{Kibler}. This reference and references therein are a neraly complete background for the mathematical problem and its physical applications. 
%
%
%

Using a similar argument as in the proof of Theorem \ref{Insta_Per} for the Peregrine case, one can show that Kuznetsov-Ma breathers are unstable \cite{Munoz2}.
The (formally) unstable character of Peregrine and Kuznetsov-Ma breathers was well-known in the physical and fluids literature (they arise from modulational instability), therefore the conclusions from previous results are not surprising. Even tough, in water tanks and optic fiber experiments, researchers were able to reproduce these waves \cite{Brunetti,Kliber0,Kibler}; if e.g. the initial setting or configuration is close to the exact theoretical solution. 


Floquet analysis has been recently done for the KM breather in \cite{Cuevas}. Concerning the variational structure of the KM breather, it is slightly more complicated than the one for Peregrine, but it has the same flavor.

\begin{thm}[\cite{AFM1}]\label{TH1}
Let $B=B_{KM}$ be a Kuznetsov-Ma breather \eqref{KM}. Then  
\be\label{Ec_KM}
\begin{aligned}
& B_{(4x)} + 3B_x^2 \bar B +(4 |B|^2-3) B_{xx}+ B^2 \bar B_{xx}  + 2  |B_x|^2 B    + \frac32 (|B|^2-1)^2 B - \beta^2(B_{xx} +  (|B|^2-1) B) =0.
\end{aligned}
\ee
\end{thm}

Note that the elliptic equation for the $P$ breather \eqref{Ec_P} is directly obtained by the formal limit $\beta\to 0$ in the $KM$ elliptic equation \eqref{Ec_KM}. This is concordance with the expected behavior of the $KM$ breather as $a\to \frac12^+$, see \eqref{KM}.


The variational structure of $B_{KM}$ goes as follows. Define 
\be\label{H_KM}
\mathcal H[u] := F[u] + \beta^2 E[u].
\ee
Here, $E$ and $F$ are given by \eqref{Mass_NLS}, and \eqref{F_NLS}, respectively. Then, for any $z\in H^2(\R;\Com)$, $\mathcal H'[ B_{KM}](z)=0.$ One has

\begin{thm}[Absence of spectral gap and instability of the KM breather, \cite{AFM1}]\label{TH5}
Let $B=B_{KM}$ be a Kuznetsov-Ma breather \eqref{KM}, critical point of the functional $\mathcal H$ defined in \eqref{H_KM}. Then
\be\label{espectro}
\begin{aligned}
\mathcal H'[B_{KM}]= 0,\quad \mathcal H''[B_{KM}](\partial_x B_{KM}) = 0, \quad \hbox{and}\quad \inf \sigma_c (\mathcal H''[B_{KM}]) < {} 0.%
\end{aligned}
\ee
\end{thm}

Note that classical stable solitons or solitary waves $Q$ \eqref{soliton} easily satisfy the estimate $\inf \sigma_c (\mathcal H_Q''[Q])>0$, where $\mathcal H_Q''$ is the standard quadratic form associated to the energy-mass-momentum variational characterization of $Q$ \cite{Weinstein}. Even in the cases of the mKdV breather $B_{mKdV}$ \cite{AM} or Sine-Gordon breather $B_{SG}$ \cite{AMP1}, one has the gap  $\inf \sigma_c (\mathcal H_{mKdV}''[B_{mKdV}])>0$ and $\inf \sigma_c (\mathcal H_{SG}''[B_{SG}])>0$. The $KM$ breather does not follow this property at all, giving another hint of its unstable character.

The above theorem shows that the $KM$ linearized operator $\mathcal H''$ has at least one \emph{embedded eigenvalue}. This is not true in the case of linear, 
real-valued operators with fast decaying potentials, but since  $\mathcal H''$ is a matrix operator, this is perfectly possible. 
Recall that if $\inf \sigma_c (\mathcal H''[B_{KM}]) > 0$, 
then the KM could perfectly be stable, a contradiction.

\subsection{The Akhmediev breather}\label{Sect:8}

Recall the Akhmediev breather \eqref{Ak}. Note that $B_A$ is a $\frac{2\pi}{a}$-periodic in space, localized in time smooth solution to \eqref{NLS0}, with some particular properties at spatial infinity. In the limiting case $a \uparrow \frac12$ one can recover the Peregrine soliton \eqref{P}. Moreover, one has
\be\label{Asymptotic}
\lim_{t\to\pm \infty} \|B_A(t,x) - e^{\pm i \theta}e^{it}\|_{H^1_\sharp}=0, \qquad e^{i\theta}=1-\al^2 -i\beta .
\ee
The instability of $B_A$ goes as follows. Once again, being $B_A$ unstable, it does not mean that it has no structure at all.

\begin{thm}[\cite{AFM2}]\label{Insta_A}
The Akhmediev breather \eqref{Ak} is unstable under small perturbations in $H^s_\sharp$, $s>\frac12$. Also, it is a critical point of the functional $\mathcal H[u]:= F_A[u] - \al^2 E_A[u]$, i.e.  $\mathcal H'[B_A][w]=0$ for all $w\in H^2_\sharp$. In particular, for each $t\in\R$, $B_A(t,x)$ satisfies the nonlinear ODE
\be\label{Ec_A}
 A_{(4x)} + 3A_x^2 \bar A +(4 |A|^2-3) A_{xx}+ A^2 \bar A_{xx}  + 2  |A_x|^2 A   + \frac32 (|A|^2-1)^2 A + \alpha^2(A_{xx} +  (|A|^2-1) A) =0.
\ee
\end{thm}

The proof of Theorem \ref{Insta_A} uses \eqref{Asymptotic} in a crucial way: a modified Stokes wave is an attractor of the dynamics around the Akhmediev breather for large time. 

\begin{rem}\label{rem:twosoliton}
We finally remark that the three above discussed breathers, $B_P, B_{KM}, B_A$ are all related 
to the two-soliton solutions of NLS on the plane-wave background (see Chapter 3 in \cite{AC}). Indeed, 
the stationary Lax-Novikov equations for all breathers belong to the same family  \eqref{Ec_P}, \eqref{Ec_KM}, 
\eqref{Ec_A} (see the recent paper \cite{CPW} for a detailed discussion about the stationary Lax-Novikov equations).
Similar conclusions are expressed in \cite{AFM1}.
\end{rem}

\section{Conclusions}\label{Sect:7}

We have reviewed stability properties of three NLS solutions with nonzero background: Peregrine, Kuznetsov-Ma and Akhmediev breathers. Working in associated energy spaces, with no additional decay condition, this review also characterizes spectral properties of each of them. According to the definition of stability, no NLS \eqref{NLS0} breather seems to be stable, not even in larger spaces. The instability is easily obtained from the fact that each breather converges on the whole line, as time tends to infinity, towards the Stokes wave. If the solutions were stable, this would imply that each breather is the Stokes wave itself. Some deeper connections between the stability of breathers and the nonzero background (modulational instability) are highly expected, but it seems that no proof of this fact is in the literature. Maybe B\"acklund transformations, in the spirit of \cite{AM1,AMP3,MP}, could help to give preliminary answers, and rigorous IST methods such as the ones in \cite{Gong1,Gong2} may help to solve this question. Finally, the dichotomy blow up/global well-posedness, and ill-posedness for large data in NLS \eqref{NLS0} with nonzero background, are interesting mathematical open problems to be treated elsewhere.

\medskip

\noindent
{\bf Acknowledgments.} We thank both referees for several comments, criticisms and the addition of several unattended references that led to an important improvement of a previous version of this manuscript.




%
%



\end{document}